\documentclass{amsart}                   
\usepackage{amsfonts}
\usepackage{amsmath}
\usepackage{amscd}
\usepackage{amssymb}
\usepackage[dvips]{graphicx}
\usepackage[dvips]{epsfig}
\usepackage{color}

\def\mf#1{\mathfrak{#1}}
\def\mb#1{\mathbb{#1}}
\def\mc#1{\mathcal{#1}}

\def\e#1{\text{\tt e}\sub{#1}}
\def\Base#1{\left\{\ \e0,\dots,\e{#1}\ \right\}}

\def\sup#1{^{^{#1}}}
\def\sub#1{_{_{#1}}}
\def\R#1{\mb{R}\sup{#1}}                      
\def\N{\mb{N}}
\def\Z{\mb{Z}}
\def\im{\text{\rm Im}}              

\def\ARROW#1{
\text{\begin{picture}(35,18)(15,15)         
            \put(30,27){$\sub{#1}$}
            \put(18,18){\vector(1,0){30}}
\end{picture}}}

\def\LARROW#1{
\text{\begin{picture}(35,18)(15,15)         
            \put(30,27){$\sub{#1}$}
            \put(48,18){\vector(-1,0){30}}
\end{picture}}}

\def\e#1{\text{\tt e}\sub{#1}}

\def\qbasic#1{\left[ #1\right]_{_q}}
\def\qfac#1{\qbasic{#1 }!\ }
\def\qcomb#1#2{\left[\begin{array}{c}  #1 \\ #2 \end{array}\right]_{_q}}

\def\morfg#1#2{\text{Hom}\left({#1},{#2}\right)}

\def\qComplex#1#2{\ \sub{q}\!C\sub{#1}\left({#2}\right)}

\def\qChains#1#2{\ \sub{q}\!SC\sub{#1}\left({#2}\right)}

\def\Hom#1#2#3{H\sub{#1}\sup{#2}\left(#3\right)}               
              
\def\qHom#1#2#3{\ \sub{q}\!H\sub{#1}\sup{#2}\left(#3\right)}          
\def\qRedHom#1#2#3{\ \sub{q}\!\widetilde{H}\sub{#1}\sup{#2}\left(#3\right)}

\newcounter{numero}

\newcounter{letra}

\newcounter{romnumero}

\newcounter{bibnumero}

\newcounter{romano}

\newtheorem{teo}{Theorem}[subsection]                  
\newtheorem{lema}[teo]{Lemma}
\newtheorem{prop}[teo]{Proposition}
\newtheorem{cor}[teo]{Corollary}

\def\bteo{\begin{teo}}
\def\eteo{\end{teo}}
\def\bprop{\begin{prop}}
\def\eprop{\end{prop}}
\def\bcor{\begin{cor}}
\def\ecor{\end{cor}}
\def\blema{\begin{lema}}
\def\elema{\end{lema}}

\theoremstyle{definition}                           
\newtheorem{definition}[teo]{Definition}
\newtheorem{definitions}[teo]{Definitions}
\newtheorem{ejems}[teo]{Examples}
\newtheorem{ejem}[teo]{Example}
\newtheorem{falsa}[teo]{}

\def\bdeff{\begin{definition}\rm }
\def\edeff{\hfill$\square$ \end{definition}}
\def\bdefs{\begin{definitions}\rm }
\def\edefs{\end{definitions}}
\def\bejem{\begin{ejem}\rm }
\def\eejem{\end{ejem}}
\def\bejems{\begin{ejems}\rm }
\def\eejems{\end{ejems}}
\def\bfalsa{\begin{falsa}\rm }
\def\efalsa{\end{falsa}}

\theoremstyle{remark}                               
\newtheorem{obs}[teo]{Remark}

\def\bobs{\begin{obs} }
\def\eobs{\end{obs}}

\newenvironment{dem}{ [{\it Proof\/}]\rm\hskip3mm }{\hfill$\square$\vskip5mm}       
\newenvironment{sketch}{ [{\it Sketch of the Proof\/}]\rm\hskip3mm }{\hfill$\square$\vskip5mm}

\def\bdem{\begin{dem}}
\def\edem{\end{dem}}

\def\bsketch{\begin{sketch}}
\def\esketch{\end{sketch}}

\begin{document}

\author{M. Angel}
\address{Grupo de \'Algebra y L\'ogica, Escuela de Matem\'aticas. Universidad Central de Venezuela.\vskip1mm
Facultad de Ciencias, Escuela de Matem\'aticas. Ciudad Universitaria, Av. Los Ilustres.\vskip1mm
Caracas 1010 Venezuela. Tlf. (+58212)6051199}
\email{mauricio.angel@ciens.ucv.ve}

\author{G. Padilla}
\address{Departamento de Matem\'aticas, Edificio 404, Ofic. 315\vskip1mm
Universidad Nacional de Colombia, Facultad de Ciencias\vskip1mm
Carrera 30, calle 45. Bogot\'a - Colombia. Tlf. (+571)3165000. Ext. 13166}
\email{gipadillal@unal.edu.co}

\title{$q$-Analog Singular Homology of Convex Spaces}
\dedicatory{\tiny Dedicated to Professor A. Reyes,\\ in the occasion of his 76th birthday.}
\date{15/02/2011}

\begin{abstract}

In this article we study some interesting properties of the $q$-Analog singular homology, which is a generalization of the usual singular homology, suitably adapted to the context of $N$-complex and amplitude homology \cite{kapranov}. We calculate the $q$-Analog singular homology of a convex space. Although it is a local matter; this is an important step in order to understand the presheaf of $q$-chains and its algebraic properties. Our result is consistent with those of Dubois-Viol\`ette \& Henneaux \cite{dubois3}. Some of these results were presented for the XVIII Congreso Colombiano de Matem\'aticas in Bucaramanga, 2011.

\end{abstract}
\maketitle

\section*{Introduction}
The fact that singular homology satisfies the homotopy axiom is a well known result of topological algebra. It can be understood in several ways. From the topological scope it asserts that any topological space that is homotopic to a single point, must have no topological holes. More than this, homotopic spaces have the same singular homology and homotopic maps induce the same maps between the respective homology groups. A customary proof can be carried out by means of these mathematical facts, 
\begin{enumerate}
	\item The cone construction \cite[p.33]{dold}.
	\item A Leibnitz rule for the convex product of singular chains \cite[p.220]{bredon2}.
	\item The double composition of border map $\partial$ vanishes, i.e. $\partial\sup2=0$, which means that singular chains constitute a usual
	chain complex. 
\end{enumerate}
On the other hand, the theory of $N$-complexes has raised in the last years as a new homology theory with a broad field of applications in quantum physics \cite{dubois3}. Let $N\geq 3$ be a prime integer. A $N$-complex is a graded module whose border map $\partial$ vanishes in the $N$-th composition, i.e. $\partial\sup{N}=0$. The $m$-amplitude homologies are defined for $1\leq m\leq N-1$; see \cite{dubois, kapranov}. For instance, take a complex $N$-th root of the identity,  $q\in\mb{C}$; i.e.  $q^{N}=1$. Then there can be defined $q$-simplicial chains, as singular chains that are linear combinations of singular simplexes where the constants are taken on the ring $\Z[q]$ and the border map is adequately adapted. Several examples will be treated here below. \vskip2mm

The main result in this article is that any convex Euclidean space has the same $q$-Analog singular homology of a singleton. This is a consequence of the algebraic structure induced by the border map and the combinatorial properties of $q$-numbers. In order to prove this,
\begin{enumerate}
	\item We use the fact that $\partial\sup{N}=0$, i.e. $q$-Analog singular chains are a graded $N$-differential module.
	\item We extend the cone construction to a convex product for the q-Analog singular homology.
	\item We obtain a $q$-Leibnitz rule for the convex product and a formula for the Newton's polynomials.
	\item We construct a geometric $N$-homotopy operator by means of the convex product.
\end{enumerate}
An open question we hope to answer in the future is to demonstrate that $q$-Analog singular chains satisfy the Mayer-Vietoris property.\vskip2mm

The article has been organized as follows. In the sections \S1,\S2  we summarize some usual facts of $q$-numbers and $N$-complexes. Section \S3 is devoted to $q$-singular chains and more examples. In section \S4 we define the convex product and show the Leibnitz rule. The last section is devoted to prove the homotopy axiom for $q$-Analog singular homology, which is our main result.

\section{$q$-numbers}
Recall the definition of $q$-numbers and some of their properties \cite{kapranov}.

\subsection{$q$-numbers}\label{def q numbers}
    Let $q\in\mb{C}$ be a complex non trivial
    $N$-th root of the identity  i.e.  $q^N=1$ and $q\neq1$. In the classical
    literature $N$ is assumed to be a prime integer and $q=\exp\left(2\pi i/N\right)$, see \cite{kapranov}. The {\bf basic $q$-numbers}  are
    \begin{equation}\label{eq basic numbers}
    			\qbasic{k}=\frac{1-q^k}{1-q}=1+q+\dots+q^{(k-1)}\hskip2cm \forall k\in\N
	\end{equation}
	Notice that $\qbasic{N}=0$.  The {\bf  $q$-factorial} numbers are
    \begin{equation}\label{eq basic factorials}
    			\qfac{k}=\qbasic{1}\cdot\qbasic{2}\cdots\qbasic{k}\hskip2cm 0\leq k\leq N-1
	\end{equation}
	Finally, the {\bf $q$-combinatorial numbers} are 
    \begin{equation}\label{eq basic combinatorials}
    			\qcomb{k}{l}=\frac{\qfac{k}}{\qfac{l}\qfac{(k-l)}}\hskip2cm \forall 0\leq l\leq k\leq N-1
	\end{equation}
Since $N\geq 2$ the polynomial $q\sup{N}-1=0$ is irreducible in $\R{}$, so $\R{}[q]=\mb{C}$ is the field of
complex numbers. In particular, $q$-numbers $\qbasic{k}\neq0$ have multiplicative inverse in $\R{}[q]$.
The following properties follow from the definition of $q$-numbers, we leave the details to the reader.
\blema\label{lema qcomb}
	Let $1\leq k\leq n,m\leq N-1$. Then,\vskip1mm
	\begin{enumerate}
		\item  $\qbasic{m+n}=\qbasic{m}+q\sup{m}\qbasic{n}$.\vskip1mm
		\item  If $n$ is prime relative to $N$, then $\qbasic{n}$ is a unit in $\Z[q]$; and 
			its multiplicative inverse is $[a]\sub{q\sup{n}}$ where $an+bN=1$ for some 
			integers $a,b$.\vskip1mm
		\item $\qcomb{n}{k}+q\sup{k+1}\qcomb{n}{k+1}=\qcomb{n+1}{k+1}=\qcomb{n}{k+1}+q\sup{n-k}\qcomb{n}{k}$.
		\item $\qfac{n}=\underset{_{\sigma\in S\sub{n}}}{\sum} q\sup{\text{\rm sgn}(\sigma)}$ where $S\sub{n}$
			is the $n$th symmetric group and $\sigma$ runs over all permutations of $n$ elements.
	\end{enumerate}
\elema	

\section{$N$-complexes}\label{section N-modules}
Let us fix a positive integer $N\geq2$ and a principal ideal domain $(R,+,\cdot,1)$ as the underlying ring of constants (usually we will take $R=\Z[q]$). A $N$-complex is a generalization of usual chain complexes, and presents a similar behavior taking into account the integer $N$, which is called the amplitude of the complex \cite{dubois,kapranov}. 

\subsection{$N$-complexes}\label{def Ndiffmodule}
       A {\bf $N$-complex} is a pair $(M,\partial)$ such that $M$ is a module and $M\ARROW{\partial}M$ 
	is a linear endomorphism such that the $N$-th composition $\partial\sup{N}=0$ vanishes. We call $\partial$ the {\bf border map}. For any integer $1	
	\leq m    \leq N-1$, we consider the submodules 
    \[
    			M\ARROW{\hskip-4mm \partial\sup{(N-m)}}M\ARROW{\partial\sup{m}}M\hskip2cm B\sub{m}(M)=\im\left(\partial\sup{(N-m)}\right)
			\subset \ker\left(\partial\sup{m}\right)=Z\sub{m}(M)
    \]
    An element of $Z\sub{m}(M)$ (resp. $B\sub{m}(M)$) is a $m$-{\bf amplitude} {\bf cycle} (resp.  {\bf border}) .   
    {\bf The homology of $M$ with amplitude $m$} is the quotient module 
    \[	
    		\Hom{m}{}{M}=\frac{Z\sub{m}(M)}{B\sub{m}(M)}
	\] 
    The {\bf total homology} of $M$ is the graded module 
    \[
    		\Hom{}{}{M}=\{\Hom{m}{}{M}: 1\leq m \leq N-1\}
	\]
	A {\bf morphism of $N$-complexes} $(M,\partial)\ARROW{f}(M',\partial')$ is a linear morphism $f$  such that
	$f\partial =\partial'f$. The induced arrow is well defined on each amplitude homology $\Hom{m}{}{M}\ARROW{f}\Hom{m}{}{M'}$, and 
	passes to the total homology $\Hom{}{}{M}\ARROW{f}\Hom{}{}{M'}$. \vskip2mm
	
	For any short exact sequence of $N$-complexes
	\[
			0\ARROW{}M\ARROW{\alpha}M'\ARROW{\beta}M''\ARROW{}0
	\]
	there is a version of the snake lemma, and a connecting morphism $\Hom{m}{}{M}\ARROW{\partial}\Hom{N-m}{}{M}$
	from which arises an exact hexagon,
		   \begin{center}
    \begin{picture}(30,90)(15,15)         
            \put(-150,50){$\Hom{m}{}{M}$}
            \put(-50,80){$\Hom{m}{}{M'}$}
            \put(50,80){$\Hom{m}{}{M''}$}
            \put(150,50){$\Hom{N-m}{}{M}$}
            \put(50,20){$\Hom{N-m}{}{M'}$}
            \put(-50,20){$\Hom{N-m}{}{M''}$}

            \put(-115,62){{\vector(3,1){60}}}
            \put(-5,85){{\vector(1,0){50}}}
            \put(90,82){{\vector(3,-1){60}}}
            \put(165,43){{\vector(-3,-1){60}}}
            \put(48,23){{\vector(-1,0){40}}}
            \put(-55,25){{\vector(-3,1){60}}}

            \put(-100,75){$^{\alpha}$}
            \put(20,90){$^{\beta}$}   
            \put(130,75){$^{\partial}$}            
            \put(130,20){$^{\alpha}$}            
            \put(20,5){$^{\beta}$}                     
            \put(-110,20){$^{\partial}$} 
                        
            \end{picture}
    \end{center}\vskip2mm

\subsection{Graded $N$-complexes}\label{Ncomplexes}
    A {\bf graded $N$-differential
    module} is a pair  $\left(M\sub{*},\partial\right)$ such that \linebreak $M\sub{*}=\{M\sub{k}:k\in\Z\}$ is a graded module  and  $\partial$ is a $(-1)$-graded linear
    endomorphism
    \[
    			M\sub{k}\ARROW{\partial}M\sub{k-1}\hskip2cm
			i\in\Z
	\]			
	such that $\partial\sup{N}=0$. The properties of $N$-complexes can be extended to the graded case. 
	The amplitude homology is now a     bigraded module 
    \[
			\Hom{}{}{M}=\left\{\Hom{m,k}{}{M}:1\leq m\leq N-1,\ k\in\Z\right\}
    \]
    depending on the amplitude $m$ and the degree $k$. The inclusion $i$ and the border map $\partial$ induce, respectively,
    well defined maps in the bigraded homology. 

\subsection{Examples}\label{examples}
	\begin{enumerate}
		\item Any finite sequence of modules and morphisms
		\[
				0\LARROW{}M\sub{1}\LARROW{_{\partial\sub1}}M\sub{2}\LARROW{_{\partial\sub2}}\cdots
				\LARROW{_{\hskip-2mm \partial\sub{N-2}}}
				M\sub{N-1}
				\LARROW{_{\hskip-2mm \partial\sub{N-1}}}M\sub{N}\LARROW{}0
		\]
		is a graded $N$-complex.
		\item With a little abuse of notation let us write 
		\[
			\Z[q]\LARROW{_{\qbasic{n}}}\Z[q]
		\] 
		for the linear function 
		that maps any element $\alpha\in\Z[q]$ to $\qbasic{n}\cdot\alpha$. According to \S\ref{lema qcomb}-(2), 		
		since $N$ is prime, $\qbasic{n}\neq0$ 
		has a multiplicative inverse in $\Z[q]$ for $1\leq n\leq N-1$. The above map is a module 
		isomorphism between free $\Z[q]$-modules.
		{\small\[
				0\LARROW{}\Z[q]\LARROW{_{\qbasic{2}}}\Z[q]\LARROW{_{\qbasic{3}}}\Z[q]\LARROW{}\cdots
				\LARROW{_{\hskip-2mm\qbasic{N-1}}}\Z[q]
				\LARROW{0}0\LARROW{}0\LARROW{}\cdots
		\]} 
		is a $N$-complex; we use to denote it by $\left(\Z[q], \qbasic{*}\right)$. 
		A straightforward calculation shows that 
		\[
			\Hom{m,n}{}{\Z[q],\qbasic{*}}=\left\{
			\begin{array}{lll}
				\Z[q] & & 1\leq n=m\leq N-2 \\
				0     & & \text{else}
			\end{array}\right.
		\]
		\item One can construct $N$-differential modules with smooth differential forms on $\R{n}$; see	
		 \cite{dubois, kapranov}. There are also  $N$-complexes
		with geometric singular chains on any topological space. For more details see the next sections.
	\end{enumerate}

\subsection{Homotopy of $N$-complexes}\label{Nhomotopy}
    Given any two morphisms of $N$-differential modules
    $M\ARROW{\hskip-2mm f,g}M'$, we say that they are {\bf homotopic} and write $f\sim g$ iff there is
    a sequence of morphisms of modules $M\ARROW{K\sub{m}}M'$, for $0\leq m \leq N-1$, satisfying 
    \begin{equation}\label{eq Nhomotopy}
    			\underset{_{m=0}}{\overset{_{N-1}}{\sum\ }} (\partial')\sup{m}K\sub{m} \partial\sup{N-m-1}  = (f-g)
    \end{equation}
  	The sequence of morphisms $K=\{K\sub{m}\}\sub{m}$ is a {\bf homotopy} from $g$ to $f$.  
    The existence of homotopies is an equivalence relation between morphisms of $N$-complexes; 
  	homotopic morphisms induce the same maps in the amplitude homologies. An alternative way to see that 
  	this is  suitable definition of homotopy between morphisms of differential $N$-modules is to 	
  	to follow \cite{tanre}[p.4-5]. Consider, for any
  	pair of $N$-differential graded modules $(M,\partial)$ and $(N,\delta)$, the graded module
  	$\morfg{M}{N}$ with the $N$-differential operator given by  
		\begin{equation}\label{eq Dif on Hom(M,N)}
			\mf{D}(f)=\underset{_{i=0}}{\overset{_{N-1}}{\sum\ }} q\sup{i(\deg(f)+1)}\delta\sup{i}f\partial\sup{N-i-1}
		\end{equation}
		A morphism $M\ARROW{f}N$ is {\bf compatible with the differentials} iff it is a $\mf{D}$-cycle, and then it 
		induces a well defined morphism on the $k$-amplitude homologies
		$\Hom{k}{}{M}\ARROW{f}\Hom{k}{}{N}$ for $1\leq k\leq N-1$. Then, two differential morphisms $f,g$ (with $\deg(f)=\deg(g)=0$ 
		as above) are {\bf homotopic} iff 	
		their difference $f-g$ is a $\mf{D}$-border in $\morfg{M}{N}$. This happens iff there exists a morphism 
		$M\ARROW{K}N$ such that $\deg(k)=(N-1)$ and $(f-g)=\mf{D}(k)$. Notice that the morphism $K$ has degree $\deg(K)=N-1$.

\section{$q$-Chains}\label{section q-homology}

\subsection{The $N$-complex of $q$-chains on a simplicial set}\label{subsection Nq-chains}
  Recall the construction of simplicial $q$-chains \cite{dubois,kapranov}.
A {\bf simplicial set} is a family  of non-empty sets and maps
\[
	X\sub{n+1}\ARROW{\partial\sub{i}}X\sub{n}\hskip2cm  0\leq i\leq n, \ n\in\N
\]
such that their compositions (\footnote{We write $fg$ for the composition $f(g(x))$ on each $x$ where it makes sense.}) satisfy
\[
	\partial\sub{i}\partial\sub{j}=\partial\sub{j}\partial\sub{i+1}\hskip2cm\forall j\leq i
\]
An element of $X\sub{n}$ is a basic chain of dimension $n$. Let $N$ and $q$ be as in \S\ref{def q numbers}.  Take the polynomial extension $\Z[q]$ as the ring of constants. The {\bf $(N,q)$-complex generated  by $X$} is the graded free $\Z[q]$-module that on each degree $n$ is spanned by $X\sub{n}$ as a linear basis.
\[
	\qComplex{n}{X}=\Z[q]\left\langle X\sub{n}\right\rangle =\underset{_{x\in X\sub{n}}}{\oplus} \Z[q]\cdot x\hskip2cm
	n\in\N
\]
As usual we assume the convention $\qComplex{n}{X}=0$ for $n<0$. The {\bf border map} is the graded linear morphism
\[
	\qComplex{n}{X}\ARROW{\partial}\qComplex{n-1}{X}\hskip1cm
	\partial=\underset{_{i=0}}{\overset{_{n}}{\sum\ }}q\sup{i}\partial\sub{i}
\]
We must check that our definition makes sense.

\blema\label{lema iteration rule border map}
	{\bf [Iteration rule for the border map]} The following equality holds
	\[
		\partial\sup{k}=\qfac{k}\cdot\underset{_{i\sub1\leq\cdots\leq i\sub{k}}}{\ \sum\ }
		q\sup{i\sub1+\cdots+i\sub{k}}
		\partial\sub{i_k}\cdots\partial\sub{i_1};
		\hskip1cm
		0\leq k\leq N
	\]
	Therefore, $\left(\qComplex{*}{X},\partial\right)$ is a graded $N$-complex.
\elema
\bdem
	Apply the definition of the border map $\partial$ and property \S\ref{lema qcomb}-(4).
	See \cite{kapranov}.
\edem
In particular, since $\qbasic{N}=0$ we get $\partial\sup{N}=0$, so $\qComplex{*}{X}$ is a $N$-complex.

\subsection{Singular $q$-chains}\label{def singular q-chains}
	 A geometric realization is given by the  $N$-complex of {\bf Singular $q$-chains}. 
	For each integer $n\in\N$ we write $\Delta\sup{n}$
	 for the standard
	$n$-simplex, i. e. the convex hull generated on $\R{n+1}$ with the standard basis $\Base{n}$.
	A linear map $\Delta\sup{n}\ARROW{L}\Delta\sup{m}$ is determined by its values on $\e0,\dots,\e{n}$;
	we write $L=\langle x\sub0,\dots,x\sub{n}\rangle$ to mean that
	$x\sub{i}=L\left(\e{i}\right)$ for $i=0,\dots,n$. Take
	\[
			\Delta\sup{n}\ARROW{\lambda\sub{j}}\Delta\sup{n+1}\hskip1cm
			\lambda\sub{i}=\langle \e0,\dots,\widehat{\e{j}},\dots,\e{n+1}\rangle
			\hskip1cm j=0,\dots,n
	\]\\ 
	where $\widehat{\e{j}}$ means to omit the element $\e{j}$. Given a topological space
	$X\neq\emptyset$ we define $X\sub{n}$ as the set of all
	continuous maps $\Delta\sup{n}\ARROW{\sigma}X$. An element of $X\sub{n}$
	is a {\bf simplex on $X$}. For each $0\leq j\leq n$ the {\bf $j$th-face map} 
	$X\sub{n}\ARROW{\partial\sub{j}}X\sub{n-1}$ is given by the composition
	$\partial\sub{j}(\sigma)=\sigma\lambda\sub{j}$. This family is a simplicial set
	in our previous sense. The {\bf $N$-complex of singular $q$-chains} on a topological space $X$ 
	\[
		\qChains{n}{X}\ =\ \qComplex{n}{\left\{X\sub{n}:n\in\N\right\}\cup
		\left\{\partial\sub{j}:X\sub{n}\ARROW{}X\sub{n-1} : 0\leq j\leq n,\ n\in \N\right\}}
	\]  
	is the $(N,q)$-complex generated by the singular $q$-simplexes and face maps. An element 
	$\xi\in\qChains{n}{X}$ is a {\bf singular $q$-chain of dimension $n$}; it can be written  a
	linear combination $\xi=a\sub1\sigma\sub1+\cdots+a\sub{r}\sigma\sub{r}$ where
	each $a\sub{i}\in Z[q]$ is a polynomial and each $\sigma\sub{i}\in X\sub{n}$
	is a simplex of dimension $n$ on $X$. We also write $n=\dim(\xi)$.
	The standard {\bf singular $(N,q)$-homology} of $X$ is the homology of this $N$-complex 
	\[
		\qHom{m,n}{}{X}\ =\ \qHom{m}{}{\qChains{n}{X}}\hskip2cm \ 1\leq m\leq N-1,\ n\in\N
	\]

\subsection{Example: $q$-homology of a point}\label{example point homology}
	If $P=\{p\}$ is a single point; then $P\sub{n}=\{\sigma\sub{n}\}$ where $\Delta\sup{n}\ARROW{\sigma\sub{n}}P$ is 
	the constant map. The module 
	\[
		\qChains{n}{P}=\Z[q]\cdot\sigma\sub{n}\cong\Z[q]
	\]
	is isomorphic to the ring of constants $\Z[q]$ through the change of basis $\sigma\sub{n}\mapsto 1$. 
	All face maps $\partial\sub0=\cdots=\partial\sub{n}$ 
	coincide. The border operator $\qChains{n}{P}\ARROW{\partial}\qChains{n-1}{P}$ 
	is the zero map for $n=0$. For $n\geq 1$ 
	\[
		\partial\left(a\sigma\sub{n}\right) =\left(\partial\sub0+q\partial\sub1+q\sup2\partial\sub2+\cdots+q\sup{n}\partial\sub{n}\right)
		\left(a\sigma\sub{n}\right)=
		\left(1+\cdots+q\sup{n}\right)\partial\sub0\left(a\sigma\sub{n}\right)
		=\qbasic{n+1}a\sigma\sub{n-1}
	\]
	 can be seen as the multiplication by the element $\qbasic{n+1}$;
	\[
		\Z[q]\ARROW{\partial}\Z[q]\hskip1cm
		\partial(a)=\qbasic{n+1}\cdot a
	\]
	 It vanishes when $n+1$ a positive multiple of $N$.  In any other case  $\qbasic{n+1}\neq0$ is
	 a unit in $\Z[q]$; see \S\ref{lema qcomb}-(2), so $\partial$ is a module isomorphism (though not a ring isomorphism).  Therefore,
	\begin{equation}\label{eq point homology}
		\qHom{m,n}{}{P}=\left\{
		\begin{array}{cl}
						\Z[q] &   0\leq n=m-1\leq N-2 \\
						& \\
					0 & \text{else}\\
		\end{array}
		\right.
	\end{equation}
	coincides with the amplitude homology of the $N$-complex given in
	the first examples \S\ref{examples}-(2).

\subsection{Exact sequence of a pair}\label{Exact sequence of a pair}
	 Given a topological space $X$ and a subspace $A\subset X$; we consider as usual the short exact sequence
	 \[
	 	0\ARROW{}\qChains{n}{A}\ARROW{}\qChains{n}{X}\ARROW{}\qChains{n}{X,A}\ARROW{}0
	 \]
	 The exact hexagon of \S\ref{def Ndiffmodule}-(6) splits to a long exact sequence
	{\small	   \begin{center}
    \begin{picture}(30,90)(15,15)         
            \put(-175,87){$\cdots$}
            \put(-115,85){$\qHom{m,n}{}{A}$}
            \put(-10,85){$\qHom{m,n}{}{X}$}
            \put(95,85){$\qHom{m,n}{}{X,A}$}
            \put(180,50){$\qHom{N-m,n-m}{}{A}$}
            \put(95,20){$\qHom{N-m,n-m}{}{X}$}
            \put(-40,20){$\qHom{N-m,n-m}{}{X,A}$}
            \put(-150,20){$\qHom{m,n-2m}{}{A}$}
            \put(-190,22){$\cdots$}

            \put(-160,90){{\vector(1,0){45}}}
            \put(-60,90){{\vector(1,0){45}}}
            \put(45,90){{\vector(1,0){50}}}
            \put(140,82){{\vector(3,-1){60}}}
            \put(220,43){{\vector(-3,-1){50}}}
            \put(90,23){{\vector(-1,0){45}}}
            \put(-42,25){{\vector(-1,0){45}}}
            \put(-145,25){{\vector(-1,0){30}}}

            \put(-50,90){$^{r}$}
            \put(65,90){$^{\delta}$}   
            \put(180,75){$^{\partial}$}            
            \put(170,35){$^{r}$}            
            \put(70,25){$^{\delta}$}                     
            \put(-70,25){$^{\partial}$} 
                        
            \end{picture}
    \end{center}}\vskip2mm 
    In the sequel, given an exact sequence from a splitted hexagon as above we will just write
    \begin{equation}\label{eq exact seq of the pair}
    		\cdots\ARROW{}\qHom{m,n}{}{A}\ARROW{}\qHom{m,n}{}{X}\ARROW{}\qHom{m,n}{}{X,A}\ARROW{\partial}
		\qHom{N-m,n-m}{}{A}\ARROW{}\cdots
    \end{equation}
    for short. In particular, this one is the {\bf $(N,q)$-homology sequence of the pair $(X,A)$}. 
    There is also a $(N,q)$-homology sequence of a triple $(X,A,B)$
    \begin{equation}\label{eq exact seq of the pair}
    		\cdots\ARROW{}\qHom{m,n}{}{A,B}\ARROW{}\qHom{m,n}{}{X,B}\ARROW{}\qHom{m,n}{}{X,A}\ARROW{\partial}
		\qHom{N-m,n-m}{}{A,B}\ARROW{}\cdots
    \end{equation}
    As usual, the connecting morphism is obtained by chasing
    in the diagram.	
    
\section{Convex product}
	Now we extend the usual cone construction \cite[p. 38]{dold} to a convex product,
	this will be the operation between $q$-chains in order to have a geometric 
	$N$-homotopy. Our goal is to construct a homotopy operator $K$ as in \S\ref{Nhomotopy}
	from the index map to the identity map in $\R{N-1}$. Since the cone constructions increases 
	the dimension in 1, a first attempt should be to iterate the conification from $N-2$ different
	affinely independent chosen points. An easier way is to take a convex combination between
	any two different singular simplexes, we develop this idea.
	
\subsection{Convex product}\label{subsection convex product}
	Suppose that $X\subset\R{d}$ is a convex subspace. Given two simplexes 
	\[
			\Delta\sup{m}\ARROW{\tau}X \LARROW{\sigma}\Delta\sup{n}
	\] 
	and a point
	\[
			(\alpha;\beta)=(\alpha\sub0,\dots,\alpha\sub{m};\beta\sub{0},\dots,\beta\sub{n})
		\in\Delta\sup{m+n+1}
	\] 
	write  $|\alpha|=\alpha\sub0+\cdots+\alpha\sub{m}$
	and $|\beta|=\beta\sub0+\cdots+\beta\sub{n}$; so $|\alpha|+|\beta|=1$.
	Consider
	\[
		\tau*\sigma:\Delta\sup{m+n+1}\ARROW{}X\hskip2cm
		\tau*\sigma(\alpha;\beta)=
		\left\{\begin{array}{lll}
			\tau(\alpha) & & |\beta|=0 \\
			&& \\
			\sigma(\beta) & & |\alpha|=0 \\
			&& \\
			|\alpha|\cdot\tau\left(\frac{\alpha}{|\alpha|}\right)+ |\beta|\cdot\sigma
			\left(\frac{\beta}{|\beta|}\right) && \text{else}
		\end{array}\right.
	\]

	The simplex $\tau*\sigma$ above is unique for each pair $(\tau,\sigma)$ so the map can be extended to a bilinear operation 
	\[
		\qChains{m}{X}\times\qChains{n}{X}\ARROW{*}\qChains{m+n+1}{X}
	\]
	For $m=0$, $\tau(\text{\tt e}\sub0)=P$ is a single point and $\tau*\sigma =P(\sigma)$ is the conification
	of $\sigma$ to the vertex $P$. In general, $\tau*\sigma$ can be thought as a convex combination of $\tau$
	and $\sigma$. The convex product satisfies nice properties with respect to the border map.

\blema\label{lema leibintz on chains}\label{prop leibintz on chains}
	{\bf [Leibnitz rule]} Let $\tau\in\qChains{m}{X}$ and $\sigma\in\qChains{n}{X}$. 
	If $mn>0$ then
	\[
		\partial(\tau*\sigma)=\partial(\tau)*\sigma + q\sup{m+1}\tau*\partial(\sigma)
	\]
\elema
\bdem
		By the bilinearity of the border map we can suppose that $\tau,\sigma$ are
		singular simplexes.
		Apply the definition of the border map, see \S\ref{subsection Nq-chains}. 
		The face maps behave as follows,
	\begin{equation}\label{eq face maps and cone prod}
		\partial\sub{i}(\tau*\sigma)=\left\{
		\begin{array}{lll}
			\left(\partial\sub{i}\tau\right)*\sigma & & 0\leq i\leq m \\
			&& \\
			\tau*\left(\partial\sub{i-m-1}\sigma\right) && m+1\leq i
		\end{array}		
		\right.
	\end{equation}
\edem

\subsection{Newton's terms}\label{Newton's truncation}
Our main goal on this \S\ is to prove a general formula for $\partial\sup{k}\left(\tau*\sigma\right)$. If 
$\tau,\sigma$ are 0-dimensional singular simplexes then, by definition of the border map; 
the border of the 1-simplex $\tau*\sigma=[\tau(\e0),\sigma(\e0)]$ is
\begin{equation}\label{border-cone of zero chains}
	\partial(\tau*\sigma)=\sigma(\e0)+q\tau(\e0)=\sigma+q\tau
\end{equation}
This is the simplest counter-example of the Leibnitz rule since, at  \S\ref{prop leibintz on chains}, the right side of the equation
vanishes. Also, if $m=\dim(\tau)=0$ and $n=\dim(\sigma)>0$, applying the definition of the border map we get
\begin{equation}\label{border-cone of zero chains 2}
	\partial(\tau*\sigma)=\sigma+q\tau*\partial(\sigma)
\end{equation}
This is exactly what happens in the usual case for $N=2$ and $q=-1$, see \cite[p.35 eq.(4.9)]{dold}; we will use this
in the sequel. Broadly speaking, since $\dim(\tau*\sigma)=m+n+1$, for $k\geq m+n+2$ all terms in a Newton's polynomial should vanish.
One can conjecture that, for $\min\{m,n\}\leq k\leq m+n+1$ some of the terms vanish and others perhaps not. 
Given a singular simplex $\Delta\sup{m}\ARROW{\tau}X$, the {\bf Newton's terms} of $\tau$ is
\[
			\mc{N}\sup{i}(\tau)=\left\{
			\begin{array}{lll}
					\partial\sup{i}(\tau) && i\leq m \\
					\qfac{m+1} && i=m+1 \\
					0 && i\geq m+2
			\end{array}
			\right.
\]
We will show the following statement,
\bprop\label{Newton truncates Leibintz}
	{\bf [Newton's polynomial]} Let $\tau\in\qChains{m}{X}$ and $\sigma\in\qChains{n}{X}$. 
	Then
	\[
		\partial\sup{k}(\tau*\sigma)=
		\underset{_{i=0}}{\overset{_{k}}{\sum}}\ q\sup{i(m+1-k+i)}\cdot\qcomb{k}{i}
				\mc{N}\sup{k-i}(\tau)*\mc{N}\sup{i}(\sigma)			
			\hskip2cm
			\forall\ k\geq 0
	\]
\eprop
In order to do so, our plan is to check the stationary behavior which begins as soon as $k>\min\{m,n\}$; 
this will be listed in a sort of lemmas called the {\bf tail formul\ae\ }. We now carry out the plan.

\blema\label{lema first tail formula}
		{\bf [Tail formula \#1]}
		Given a simplex $\Delta\sup{m}\ARROW{\tau}X$ let $T\sub{j}=\tau(\e{j})$ for $j=0,\dots,m$.
		Then
		\[
				\partial\sup{m}(\tau)=\qfac{m}\cdot
				\underset{_{j=0}}{\overset{^{m}}{\sum}} q\sup{j}T\sub{m-j}
		\]  
\elema
\bdem
	Let us apply twice the iteration rule of the border map \S\ref{lema iteration rule border map} and
	reorder the indexes. We get
	{\small
	\[
			\partial\sup2(\tau)=\underset{_{i=0}}{\overset{_{m-1}}{\sum}}
			\underset{_{j=0}}{\overset{_{m}}{\sum}} q\sup{i+j}\partial\sub{i}\partial\sub{j}(\tau)
			=(1+q)\underset{_{0\leq i\leq j\leq m-1}}{\sum} q\sup{i+j}\partial\sub{i}\partial\sub{j}(\tau)
	\]}
	Notice that the indexes $i,j$ run over $0\leq i\leq j\leq m-1$. For $\partial\sup{k}(\tau)$
	the indexes $i\sub1,\dots,i\sub{k}$ run over $0\leq i\sub1\leq\cdots\leq i\sub{k}\leq m-k+1$.
	Finally, for $\partial\sup{m}(\tau)$ the indexes $i\sub{1},\dots,i\sub{m}$ run over
	$0\leq i\sub1\leq\cdots\leq i\sub{m}=m-m+1=1$. Therefore,
	{\small
	\[
			\begin{array}{lll}
			\partial\sup{m}(\tau)
			&& =\qfac{m}\cdot\underset{_{0\leq i\sub1\leq\cdots\leq i\sub{m}\leq 1}}{\ \sum\ \ }
		q\sup{i\sub1+\cdots+i\sub{m}}
		\partial\sub{i_m}\cdots\partial\sub{i_1}(\tau) \\
		&& \\
		&& =\qfac{m}\left(
		q\sup{0+\cdots+0}\partial\sub{0}\sup{m}(\tau) + 
		q\sup{0+\cdots+0+1}\partial\sub{1}\partial\sup{(m-1)}\sub0(\tau) + \cdots +
		q\sup{1+\cdots+1}\partial\sub1\sup{m}(\tau)\right) \\
		&& \\
		&& =\qfac{m}\left(
		q\sup0\partial\sub{0}\sup{m}(\tau) + 
		q\sup1\partial\sub{1}\partial\sup{(m-1)}\sub0(\tau) + \cdots +
		q\sup{m}\partial\sub1\sup{m}(\tau)\right) \\
		&& \\
		&& =\qfac{m}\cdot
		\underset{_{j=0}}{\overset{^{m}}{\sum}} q\sup{j}T\sub{m-j}
		\end{array}
	\]}
\edem

\blema\label{lema second tail formula}
		{\bf [Tail formula \#2]}
		Given two chains $\tau\in\qChains{m}{X}$ and $\sigma\in\qChains{n}{X}$,
		\[
				\partial\left(\partial\sup{m}(\tau)*\partial\sup{n}(\sigma)\right)=
				\qfac{m+1}\partial\sup{n}(\sigma)+q\qfac{n+1}\partial\sup{m}(\tau)
				\]  
\elema
\bdem 
		By the bilinearity of the convex product and the linearity of the border map, 
		we can assume that $\Delta\sup{m}\ARROW{\tau}X\LARROW{\sigma}\Delta\sup{n}$
		are two simplexes.	 By lemma \S\ref{lema first tail formula}, let us write
		\[
				\partial\sup{m}(\tau)=\qfac{m}\underset{_{j=0}}{\overset{^{m}}{\sum}} q\sup{j}T\sub{m-j}
				\hskip2cm
					\partial\sup{n}(\sigma)=\qfac{n}\underset{_{i=0}}{\overset{^{n}}{\sum}} q\sup{i}S\sub{n-i}		\]
		in a suitable form. Then
		{\small
		\[
			\begin{array}{lll}
				\partial\left(\partial\sup{m}(\tau)*\partial\sup{n}(\sigma)\right) & 
				= \partial\left( \left( \qfac{m}\underset{_{j=0}}{\overset{^{m}}{\sum}} q\sup{j}T\sub{m-j}\right) * 
				\left( \qfac{n}\underset{_{i=0}}{\overset{^{n}}{\sum}} q\sup{i}S\sub{n-i}\right)\right) & \\
				&& \\
				& = \qfac{m}\qfac{n}\underset{_{j=0}}{\overset{^{m}}{\sum}}\  \underset{_{i=0}}{\overset{^{n}}{\sum}} q\sup{i} q\sup{j}
				\partial\left(T\sub{m-j} * S\sub{n-i}\right) \\
				&& \\
				& = \qfac{m}\qfac{n}\underset{_{j=0}}{\overset{^{m}}{\sum}}\  \underset{_{i=0}}{\overset{^{n}}{\sum}} q\sup{i} q\sup{j}
				\left(S\sub{n-i}+ q\cdot T\sub{m-j}\right) \\
				&& \\
				& = \qfac{m}\left(\underset{_{j=0}}{\overset{^{m}}{\sum}} q\sup{j}\right)\left( \qfac{n}\underset{_{i=0}}{\overset{^{n}}{\sum}} q\sup{i}
				S\sub{n-i} \right)
				+ q\cdot\qfac{n}\left(\underset{_{i=0}}{\overset{^{n}}{\sum}} q\sup{i}\right)
				    \left(\qfac{m} \underset{_{j=0}}{\overset{^{m}}{\sum}}\   q\sup{j}T\sub{m-j} \right) \\
				&& \\
				& =\qfac{m+1}\partial\sup{n}(\sigma)+q\qfac{n+1}\partial\sup{m}(\tau) & \\
			\end{array}
		\]}
		as desired.
\edem

\blema\label{lema product under border}\label{lema tail formula 1}
	{\bf [Tail formul\ae\ \#3]} Let $\tau\in\qChains{m}{X}$ and $\sigma\in\qChains{n}{X}$. 
	If $mn>0$ then\vskip3mm	
	\begin{enumerate}
		\item[] $	
		\partial\sup{k}(\partial\sup{m}(\tau)*\sigma)= \left\{
		\begin{array}{lll}
			\qfac{m+1}\qbasic{k}\partial\sup{k-1}(\sigma) 
				+ q\sup{k}\partial\sup{m}(\tau)*\partial\sup{k}(\sigma) 
			&& 1\leq k\leq n\\
			&& \\
				\qfac{m+1}\qbasic{n+1}\partial\sup{n}(\sigma) 
			+ \qfac{n+1}q\sup{n+1}\partial\sup{m}(\tau)			
			&& k=n+1\\
			&& \\
			0 && \text{\rm else}
		\end{array}\right.$\vskip3mm
	\end{enumerate}
\elema
\bdem 
	By the bilinearity of the convex product and linearity of the border map, it is enough to show it on the
	generators. Assume that $\sigma,\tau$ are simplexes. Let $T\sub{j}=\tau(e\sub{j})$ for $j=0,\dots,m$.
	By \S\ref{lema first tail formula} on $\partial\sup{m}(\tau)$, 
	equation (\ref{border-cone of zero chains 2}) at \S\ref{Newton's truncation}
	 on $T\sub{j}*\sigma$ for each $j$,
	the linearity of $\partial$ and the bilinearity of the cone-product;
	{\small\[
		\begin{array}{lll}
			\partial\left(\partial\sup{m}(\tau)*\sigma\right) & = \partial\left( \qfac{m}
			\left(\underset{_{j=0}}{\overset{^{m}}{\sum}} q\sup{j} 
			T\sub{m-j}\right)*\sigma\right)  & \\
			&&  \\	
			& =\qfac{m}\cdot
			\underset{_{j=0}}{\overset{^{m}}{\sum}} q\sup{j}\partial\left( 
			T\sub{m-j}*\sigma\right)  
			=\qfac{m}\cdot
			\underset{_{j=0}}{\overset{^{m}}{\sum}}q\sup{j}\left(\sigma + q T\sub{m-j}\right) & \\
			&&  \\				
			& =\qfac{m}\cdot\left(
			\underset{_{j=0}}{\overset{^{m}}{\sum}} q\sup{j}\right)\cdot\sigma  + 
			\qfac{m}\cdot q\left(
			\underset{_{j=0}}{\overset{^{m}}{\sum}} q\cdot\sup{j}T\sub{m-j}\right)*\partial(\sigma) & \\	
			&&  \\	
			& = \qfac{m+1}\cdot\sigma + q\cdot\partial\sup{m}(\tau)*\partial(\sigma) 	&
		\end{array}
	\]}
	This proves the equality for $k=1$; for $2\leq k\leq n$ apply this rule and use induction on $k$. 
	For $k=n+1$,  by direct calculations 
	{\small\[
		\begin{array}{lll}
			\partial\sup{n+1}\left(\partial\sup{m}(\tau)*\sigma\right) 
			& = \partial\left(\partial\sup{n}(\tau*\sigma)\right)& \\
			&& \\
			& = \partial\left(
			\qfac{m+1}\qbasic{n}\partial\sup{n-1}(\sigma) + q\sup{n}\partial\sup{m}(\tau)*
			\partial\sup{n}(\sigma) \right)	& \\
			&& \\
			& = 
			\qfac{m+1}\qbasic{n}\partial\sup{n}(\sigma) + q\sup{n}\partial\left(\partial\sup{m}(\tau)*\partial\sup{n}(\sigma) \right)
			& \\
		\end{array}
	\]}
	For the last term we now apply lemma \S\ref{lema second tail formula}. Then,
	{\small\[\hskip1cm
		\begin{array}{lll}
			\partial\sup{n+1}\left(\partial\sup{m}(\tau)*\sigma\right) 
			& = 
			\qfac{m+1}\qbasic{n}\partial\sup{n}(\sigma) + 
			q\sup{n}\left( \qfac{m+1}\partial\sup{n}(\sigma) + q\qfac{n+1}\partial\sup{m}(\tau)	\right)
			& \\
			&& \\
		& =\qfac{m+1}\qbasic{n+1}\partial\sup{n}(\sigma) 
		+ \qfac{n+1}q\sup{n+1}\partial\sup{m}(\tau)
		& \\
		\end{array}
	\]}
	as desired.
\edem

A similar expression can be obtained for $\partial\sup{k}(\tau*\partial\sup{n}(\sigma))$, though we will not
need it here.

\subsection{Proof of Proposition \S\ref{Newton truncates Leibintz}}
	We will proceed by double induction on $n+m$ and $k$. For $n+m=0$ we have $n=m=0$. Consider the following cases:
	$k=0$ which is trivial, $k=1$ which gives equation  (\ref{border-cone of zero chains}) at \S\ref{Newton's truncation},
	for $k\geq2$ we get $\partial\sup{k}(\tau*\sigma)=0$ by a dimension argument. This proves \S\ref{Newton truncates Leibintz}
	for $m+n=0$ and $k\geq0$. For $m+n>0$ fix some simplexes $\tau,\sigma$ with respective dimensions $m,n$.  Let us assume the inductive 
	hypotheses, i.e. that \S\ref{Newton truncates Leibintz} holds for any pair of simplexes
	$\tau',\sigma'$ with respective dimensions $m',n'$ such that $m'+n'<m+n$. For 
	$k=0$ there is nothing to prove. For $k\geq \dim(\tau*\sigma)+1=m+n+2$, by a dimension argument,
	 the left side of the Newton's polynomial at \S\ref{Newton truncates Leibintz} vanishes. Also, all the terms 
	 $\mc{N}\sup{k-i}(\tau)*\mc{N}\sup{i}(\sigma)$ in the right side vanish since, for any $i\leq k$, we have 
	 $i\leq n\ \Rightarrow\ k-i>m$ and $k-i\leq m\ \Rightarrow\ i>n$, so the statement holds. Hence we only have to
	 check \S\ref{Newton truncates Leibintz} for $1\leq k\leq m+n+1$. \vskip2mm

	For $k=1$ the statement of \S\ref{Newton truncates Leibintz}  is the Leibnitz rule
	\S\ref{prop leibintz on chains}. Notice that $m=\deg(\tau)$ so the power $q\sup{m+1}$ in the statement of \S\ref{prop leibintz on chains}
	depends on $\tau$; i.e. $\partial(\tau*\sigma)=\partial(\tau)*\sigma+q\sup{\deg(\tau)+1}\tau*\partial(\sigma)$. 
	Assume the inductive hypothesis for $k\leq \min\{m,n\}-1$. Then, by the linearity of the border map,
	{\small\[
			\partial\sup{k+1}(\tau*\sigma)  = \partial\left(\partial\sup{k}(\tau*\sigma)\right)
			 = \partial\left[ \underset{_{i=0}}{\overset{_{k}}{\sum}}q\sup{i(m+1-k+i)}\cdot\qcomb{k}{i}
		\partial\sup{k-i}(\tau)*\partial\sup{i}(\sigma)\right]
		 =  \underset{_{i=0}}{\overset{_{k}}{\sum}}q\sup{i(m+1-k+i)}\cdot\qcomb{k}{i}
		\partial\left(\partial\sup{k-i}(\tau)*\partial\sup{i}(\sigma)\right) 
	\]}		
	Since $k\leq \min\{m,n\}-1$, all the terms $\partial\left(\partial\sup{k-i}(\tau)*\partial\sup{i}(\sigma)\right)$ in the
	last sum satisfy the hypothesis of \S\ref{prop leibintz on chains}. Apply the Leibnitz rule to each of them. 
	We get  
	{\small\[
		\begin{array}{ll}
				\partial\sup{k+1}(\tau*\sigma) 
			& = \partial\sup{k+1}(\tau)*\sigma 
				+ 
		\underset{_{i=1}}{\overset{_{k}}{\sum}}q\sup{(i+1)(m-k+i+1)}\cdot
		\left(\qcomb{k}{i} +q\sup{(i+1)}\cdot\qcomb{k}{i+1}\right)
		\partial\sup{k-i}(\tau)*\partial\sup{i+1}(\sigma)
				
			+ \tau*\partial\sup{k+1}(\sigma) \\
		\end{array}		
	\]}
	By property \S\ref{lema qcomb}-(3) the sum of $q$-combinatorial numbers can be arranged, so
	\[
		\partial\sup{k+1}(\tau*\sigma) \ =\ 
		\underset{_{i=0}}{\overset{_{k+1}}{\sum}}q\sup{(i+1)(m-k+i+1)}\cdot
		\qcomb{k+1}{i+1}
		\partial\sup{k-i}(\tau)*\partial\sup{i+1}(\sigma)
	\]
	as desired. We have proved \S\ref{Newton truncates Leibintz} for $0\leq k\leq\min\{m,n\}$.\vskip5mm

	For $\min\{m,n\}+1\leq k\leq m+n+1$ consider the following cases.\vskip2mm 
	\begin{itemize}
	\item \underline{$m<n$:} 
	We check directly \S\ref{Newton truncates Leibintz} for $k=m+1\leq n$. Notice that
	{\tiny
	\[
		\hskip-1.5cm
		\begin{array}{lll}
			\partial\sup{m+1}(\tau*\sigma)  
			 & = \partial\left(\partial\sup{m}(\tau*\sigma)\right) =\partial\left[ \underset{_{i=0}}{\overset{_{m}}{\sum}}q\sup{i(1+i)}\cdot\qcomb{m}{i}
			\partial\sup{m-i}(\tau)*\partial\sup{i}(\sigma)\right] & \text{\S\ref{Newton truncates Leibintz} for $k=m$}\\
			   && \\
			 & = \partial\left(\partial\sup{m}(\tau)*\sigma\right) +\underset{_{i=1}}{\overset{_{m}}{\sum}}q\sup{i(1+i)}\cdot\qcomb{m}{i}
			\partial\left[ \partial\sup{m-i}(\tau)*\partial\sup{i}(\sigma)\right] &\text{linearity of $\partial$}\\
			   && \\
			& = \left[\qfac{m+1}\sigma + q\partial\sup{m}(\tau)*\partial(\sigma) \right] 
				+ 
				\underset{_{i=1}}{\overset{_{m}}{\sum}}q\sup{i(1+i)}\cdot\qcomb{m}{i}
			\left[ \partial\sup{m-i+1}(\tau)*\partial\sup{i}(\sigma)+q\sup{(1+i)}\partial\sup{m-i}(\tau)*\partial\sup{i+1}(\sigma)\right] 
			&\text{\S\ref{lema tail formula 1} and \S\ref{prop leibintz on chains}}\\
			&& \\			
			& = \qfac{m+1}\sigma + \underset{_{i=1}}{\overset{_{m}}{\sum}}
				q\sup{i^2}\cdot\left( \qcomb{m}{i-1} + q\sup{i}\cdot\qcomb{m}{i} \right)
			 \partial\sup{m-i}(\tau)*\partial\sup{i+1}(\sigma)  
				+ \tau*\partial\sup{m+1}(\sigma) & \text{Group similar terms}\\
				&& \\
			 & =\qfac{m+1}\sigma + \underset{_{j=1}}{\overset{_{m+1}}{\sum}}
				q\sup{i^2}\cdot\qcomb{m+1}{j} 
			 \partial\sup{m-i}(\tau)*\partial\sup{i+1}(\sigma) & \text{\S\ref{lema qcomb}-(3)}\\
				&& \\
		\end{array}		
	\]}
	This proves \S\ref{Newton truncates Leibintz} for $k=m+1$. Let us  
	assume again, by induction on $k$, that we have proved  it for any integer from $0$ to some 
	$k$ such that $m+1\leq k\leq n$. Then, by linearity of the border map and the inductive hypothesis,
	{\small\[
		\partial\sup{k+1}(\tau*\sigma)  
		= \partial\left(\partial\sup{k}(\tau*\sigma)\right) 
		=\partial\left[ \underset{_{i=0}}{\overset{_{k}}{\sum}}q\sup{i(m-k+1+i)}\cdot\qcomb{k}{i}
		\mc{N}\sup{k-i}(\tau)*\mc{N}\sup{i}(\sigma)\right] 
		=\underset{_{i=0}}{\overset{_{k}}{\sum}}q\sup{i(m-k+1+i)}\cdot\qcomb{k}{i}
		\partial\left( \mc{N}\sup{k-i}(\tau)*\mc{N}\sup{i}(\sigma)\right) 
	\]}
	 By definition of the Newton's terms at \S\ref{Newton's truncation},
	$N\sup{k-i}(\tau)$ vanishes for $k-i\geq m+2$. Take only take the terms satisfying $0\leq k-i\leq m+1$;
	i.e. $k-m-1\leq i\leq k$. We get,	
	{\tiny
	\[
		\hskip-1.5cm
		\begin{array}{ll}
			\partial\sup{k+1}(\tau*\sigma)  
			  & =\underset{_{i=k-m-1}}{\overset{_{k}}{\sum}}q\sup{i(m-k+1+i)}\cdot\qcomb{k}{i}
		\partial\left( \mc{N}\sup{k-i}(\tau)*\mc{N}\sup{i}(\sigma)\right)   \\
				& \\
			   \\
			   & = \qcomb{k}{k-m-1}\partial\left( \mc{N}\sup{m+1}(\tau)*\mc{N}\sup{k-m-1}(\sigma)\right)
			  +\underset{_{i=k-m}}{\overset{_{k}}{\sum}}q\sup{i(m-k+1+i)}\cdot\qcomb{k}{i}
			\partial\left(\partial\sup{k-i}(\tau)*\partial\sup{i}(\sigma)\right)   \\
			   \\
			   & \\
			   & =\qfac{m+1} \qcomb{k}{k-m-1}\partial\sup{k-m}(\sigma)
			+ q\sup{(k-m)}\qcomb{k}{k-m}  
			  \partial\left(\partial\sup{m}(\tau)*\partial\sup{k-m}(\sigma)\right)
			  \\
			   & \\
			   & \hskip4.5cm +\underset{_{i=k-m+1}}{\overset{_{k}}{\sum}}q\sup{i(m-k+1+i)}\cdot\qcomb{k}{i}
			\partial\left(\partial\sup{k-i}(\tau)*\partial\sup{i}(\sigma)\right)  \\
			   \\
		\end{array}		
	\]}
	In the last expression, apply the tail formula \S\ref{lema tail formula 1} to the second term, and
	the Leibnitz rule \S\ref{prop leibintz on chains} to the terms in the last sum.
	{\tiny
	\[
		\hskip-1.5cm
		\begin{array}{ll}
			\partial\sup{k+1}(\tau*\sigma)  
			  & =\qfac{m+1} \qcomb{k}{k-m-1}\partial\sup{k-m}(\sigma)
			+ q\sup{(k-m)}\qcomb{k}{k-m}  
			  \left(\qfac{m+1}\partial\sup{k-m}(\sigma)+q \partial\sup{m}(\tau)*\partial\sup{k-m+1}(\sigma)\right)
			   \\
			   & \\
			   & +\underset{_{i=k-m+1}}{\overset{_{k}}{\sum}}q\sup{i(m-k+1+i)}\cdot\qcomb{k}{i}
			\left(\partial\sup{k-i+1}(\tau)*\partial\sup{i}(\sigma) +q\sup{m-(k-i)+1}\partial\sup{k-i}(\tau)*\partial\sup{i+1}(\sigma)\right)  \\
		\end{array}		
	\]}
	Regroup similar terms. Apply property \S\ref{lema qcomb}-(3) on the $q$-combinatoric numbers;
	{\tiny
	\[
		\begin{array}{ll}
			\partial\sup{k+1}(\tau*\sigma)  
			   & =\qfac{m+1} \qcomb{k+1}{k-m}\partial\sup{k-m}(\sigma) 
			   + \underset{_{i=k-m+1}}{\overset{_{k}}{\sum}}\left(q\sup{i(m-k+i)}\qcomb{k}{i-1} +q\sup{i(m-k+1+i)}\cdot\qcomb{k}{i}
			\right) \partial\sup{k-i+1}(\tau)*\partial\sup{i}(\sigma) \\
			& \\
			& + q\sup{(k+1)(m+1)}\cdot\qcomb{k}{k}\partial\sup{1}(\tau)*\partial\sup{k}(\sigma) \\
			& \\
			& =\qcomb{k+1}{k-m}\mc{N}\sup{m+1}(\tau)*\mc{N}\sup{k-m}(\sigma) 
			+ \underset{_{i=k-m+1}}{\overset{_{k}}{\sum}}q\sup{i(m-k+i)}\left(\qcomb{k}{i-1} +q\sup{i}\cdot\qcomb{k}{i}
			\right) \partial\sup{k-i+1}(\tau)*\partial\sup{i}(\sigma) \\
			& \\
			& + q\sup{(k+1)(m+1)}\cdot\qcomb{k}{k}\partial\sup{1}(\tau)*\partial\sup{k}(\sigma) \\
			& \\
			& =\qcomb{k+1}{k-m}\mc{N}\sup{m+1}(\tau)*\mc{N}\sup{k-m}(\sigma) 
			+ \underset{_{i=k-m+1}}{\overset{_{k}}{\sum}}q\sup{i(m-k+i)}\qcomb{k+1}{i}\mc{N}\sup{k-i+1}(\tau)*\mc{N}\sup{i}(\sigma) \\
			& \\
			& + q\sup{(k+1)(m+1)}\cdot\qcomb{k+1}{k+1}\mc{N}\sup{1}(\tau)*\mc{N}\sup{k}(\sigma) \\
			& \\
			& =\underset{_{i=k-m}}{\overset{_{k+1}}{\sum}}q\sup{i(m-k+i)}\qcomb{k+1}{i}\mc{N}\sup{k-i+1}(\tau)*\mc{N}\sup{i}(\sigma) \\
			\end{array}
	\]}
	Include the vanishing terms of the form $\mc{N}\sup{k-i+1}(\tau)*\mc{N}\sup{i}(\sigma)$ for $0\leq i\leq k-m-1$.
	We obtain 
	\[
		\partial\sup{k+1}(\tau*\sigma)   
			=\underset{_{i=0}}{\overset{_{k+1}}{\sum}}\ q\sup{i(m-k+i)}\qcomb{k+1}{i}\mc{N}\sup{k-i+1}(\tau)*\mc{N}\sup{i}(\sigma)
	\] 
	This is the complete expression of the right term in \S\ref{Newton truncates Leibintz}
	for $k+1$. Thus we have proved the statement for $0\leq k\leq n+1$.  Finally, for $n+2\leq k\leq m+n+1$
	a similar argumentation can be carried out. The tail formul\ae\ must be used in both extremes of the sum.	\vskip2mm
	\item \underline{$m\geq n$:} We leave the details to the reader.
	\end{itemize}

\hfill$\square$\vskip5mm

\bfalsa\label{Augmentations}
	{\bf [Zeroth $q$-homology group, augmentation]}
	Since $\Delta\sup{0}=\{\e0\}$ is a singleton, each $0$-dimensional simplex $\sigma$ in $X$ can be identified
	to its image point $x=\sigma(\e0)\in X$. The $0$-th module of $q$ chains is then
	\[
			\qChains{0}{X}=\underset{_{\sigma\in X\sub0}}{\oplus}\Z[q]\cdot \sigma\cong\underset{_{x\in X}}{\oplus}\Z[q]\cdot x
	\]
	Consider the morphism
	\[
			\qChains{0}{X}\ARROW{\epsilon}\Z[q]\hskip2cm
			\underset{i}{\sum} \alpha\sub{i}x\sub{i}\ \mapsto\ \underset{i}{\sum} \alpha\sub{i}
	\]
	Given a $m$-simplex $\Delta\sup{m}\ARROW{\tau}X$ the element $\partial\sup{m}(\tau)$ is a 0-dimensional chain. 
	Let us write $P\sub{j}=\tau(\e{j})$ for $j=0,\dots,m$. Applying \S\ref{lema first tail formula} we get
	\[
		\epsilon\left(\partial\sup{m}(\tau)\right)=\qfac{m}\cdot
		\underset{_{j=0}}{\overset{^{m}}{\sum}} q\sup{j}=
		\qfac{m+1}
	\]
	In particular, for $m=N-1$ we get $\epsilon\left(\partial\sup{N-1}(\tau)\right)=0$ and
	\begin{equation}\label{eq augmentation 1}
			\qHom{1,0}{}{X}\ARROW{\epsilon}\Z[q]\hskip2cm
			[\tau]\mapsto\epsilon(\tau)
	\end{equation}
	is a well defined linear surjective morphism.

    The constant map $X\ARROW{}P$  induces a morphism of $N$-complexes
	\[
				\qChains{n}{X}\ARROW{\gamma}\qChains{n}{P}
	\] 
	called the {\bf augmentation}.    The {\bf reduced $q$-homology}
	\[
		\qRedHom{m,n}{}{X}=\ker\left\{\ \qHom{m,n}{}{X}\ARROW{\gamma}\qHom{m,n}{}{P}\ \right\}
	\]  
	is the kernel of the corresponding homology morphism. By equation (\ref{eq point homology}) at \S\ref{example point homology};
	\[
		\qHom{m,n}{}{X}=\left\{
			\begin{array}{lll}
					\Z[q]\oplus\qRedHom{m,n}{}{X} & & 1\leq n=m\leq N-2 \\
					& & \\
					\qRedHom{m,n}{}{X} & & \text{else}
			\end{array}
		\right.	
	\]  
	A reduced $q$-homology sequence of the pair
    \[
    		\cdots\ARROW{}\qRedHom{m,n}{}{A}\ARROW{}\qRedHom{m,n}{}{X}\ARROW{}\qHom{m,n}{}{X,A}\ARROW{\partial}
		\qRedHom{N-m,n-m}{}{A}\ARROW{}\cdots
    \]
	can also be deduced. 
\efalsa

\section{$q$-Analog Singular Homology of Convex Spaces}
We arrive to the main result of this article.

\subsection{The index map}\label{subsection index map}
	In complete analogy with the usual case ($N=2$, $q=-1$), the index map is, in general, the morphism  
	\[
		(\qChains{*}{X},\partial)\ARROW{\eta}\left(\Z[q], \qbasic{*}\right)
	\]
	that sends each $n$-simplex to $1\in\Z[q]$ in the corresponding degree, for $0\leq n\leq N-2$; and vanishes
	for $n\geq N-1$. 
	
\bteo\label{prop homotopy qchains}
	Let $X\subset\R{N-1}$ be a convex space. Then the  index map $\qChains{*}{X}\ARROW{\eta}\Z[q]$
	induces an isomorphism in $N$-homology.
\eteo
\bdem
	We follow essentially the same argumentation of \cite[p.38]{dold}. We will 
	define a map 
	\[
		\Z[q]\ARROW{\widehat{P}}\qChains{*}{X}
	\]
	The composition
	$\eta\widehat{P}=id$ must be the identity map on the $(N-1)$-complex $(\Z[q],[*])$; 
	so $\widehat{P}(1)=\nu\sub{n}$ will be a single singular $n$-simplex
	for $0\leq n\leq N-2$ and it will vanish for $n\geq N-1$. The other composition $\widehat{P}\eta$
	will be $N$-homotopic to the identity map $id$ on $\qChains{*}{X}$ in the sense of 
	\S\ref{Nhomotopy}. In order to explain better how we will pick the $\nu\sub{n}$'s
	we  will construct a homotopy operator
	\[
		\qChains{n}{X}\ARROW{K}\qChains{n-N+1}{X}
	\]
	and show how it works. We proceed by steps.\vskip2mm	
	\hskip-5mm$\bullet$ \underline{Definition of $K$:} Fix some singular $N-2$-dimensional simplex 
	$\Delta\sup{N-2}\ARROW{\imath} X$. Since $N$ is a prime integer, $\qbasic{k}$ is a unit in 
	$\Z[q]$ for $1\leq k \leq N-1$ and therefore $\qfac{N-1}$
	is also a unit; see \S\ref{lema qcomb}-(2).  We define 
	\begin{equation}\label{eq homotopy3}
		K(\sigma)=\frac{1}{\qfac{N-1}}\cdot(\imath*\sigma)
	\end{equation}
	Up to the correction by the constant, $K$ is essentially the convex product of $\imath$ and $\sigma$; and it can
	be uniquely extended to $\qChains{*}{\R{N-1}}$ by linearity.\vskip2mm
	\hskip-5mm$\bullet$ \underline{$K$ is a $N$-homotopy:} We verify  that $K$ satisfies \S\ref{Nhomotopy}. 
	Fix a singular simplex $\sigma\in\qChains{n}{X}$. By \S\ref {Newton truncates Leibintz} we have
	{\small\[\hskip-1cm
		\partial\sup{k}K\partial\sup{N-k-1}(\sigma)= \frac{1}{\qfac{N-1}}\cdot
		\partial\sup{k}\left(\imath*\left(\partial\sup{N-k-1}(\sigma)\right)\right)=	
		\frac{1}{\qfac{N-1}}\cdot	
		\underset{_{i=0}}{\overset{_{k}}{\sum}}\ q\sup{i(N-1-k+i)}\cdot\qcomb{k}{i}
		\mc{N}\sup{k-i}(\imath)*\mc{N}\sup{i}\left(\partial\sup{N-k-1}(\sigma)\right) 
	\]}
	Although $\partial\sup{j}(\sigma)$ is a chain and not a simplex, since $\partial\sup{i}\left(\partial\sup{j}(\sigma)\right)
	=\partial\sup{i+j}(\sigma)$ we will assume the following convention,
	\[
			\mc{N}\sup{i}\left(\partial\sup{j}(\sigma)\right)=\mc{N}\sup{i+j}(\sigma)=
			\left\{
			\begin{array}{lll}
					\partial\sup{i+j}(\sigma) && j\leq n-i \\[2mm]
					\qfac{n+1} && j=n-i+1 \\[2mm]
					0 & & \text{else}
			\end{array}
			\right.
	\]
	Therefore
	{\small\[
		\partial\sup{k}K\partial\sup{N-k-1}(\sigma)=	
		\frac{1}{\qfac{N-1}}\cdot	
		\underset{_{i=0}}{\overset{_{k}}{\sum}}\ q\sup{i(N-1-k+i)}\cdot\qcomb{k}{i}
		\mc{N}\sup{k-i}(\imath)*\mc{N}\sup{N-k-1+i}(\sigma) 
	\]}
	Taking sums in both sides, 
	{\small\[	
			\underset{_{k=0}}{\overset{_{N-1}}{\sum}}
			\partial\sup{k}K\partial\sup{N-k-1}(\sigma) 			
			= \frac{1}{\qfac{N-1}}\cdot\underset{_{k=0}}{\overset{_{N-1}}{\sum}}
			\	
			\underset{_{i=0}}{\overset{_{k}}{\sum}}\ q\sup{i(N-1-k+i)}\cdot\qcomb{k}{i}
			\mc{N}\sup{k-i}(\imath)*\mc{N}\sup{N-1-k+i}(\sigma) 
	\]}
	Let us reorder and group all similar terms taking $l=k-i$. We arrive to the following expression
	{\small\begin{equation}\label{eq homotopy4}
		\hskip-1cm
 			\underset{_{k=0}}{\overset{_{N-1}}{\sum}}
			\partial\sup{k}K\partial\sup{N-k-1}(\sigma) 			
			=  \frac{1}{\qfac{N-1}}\cdot\underset{_{l=0}}{\overset{_{N-1}}{\sum}}\ 
			\alpha\sub{l}
			\mc{N}\sup{l}(\imath)*\mc{N}\sup{N-l-1}(\sigma) 			 
	\end{equation}}
	Let us look for instance the following array of the coefficients $\alpha\sub{k,l}$ for $N=7$. The vertical sums
	of the entries in the table correspond to the values of $\alpha\sub{l}$. 
	{\small\[
	\begin{array}{|ccccccc|c|}
		\hline
		&&&&&&& k \\
		&&&&&&&  \\
		&&&&&& 1 & 0\\
		&&&&&&&  \\
		&&&&& 1 & q\sup6 & 1\\
		&&&&&&&  \\
		&&&& 1  & q\sup5\qbasic{2}  & q\sup{5} & 2\\
		&&&&&&&  \\
		&&& 1  &  q\sup4\qbasic{3}  & q\sup{3}\qbasic{3}  & q\sup{4} &  3\\
		&&&&&&&  \\
		&& 1  &  q\sup3\qbasic{4}  & q\qcomb{4}{2}  & q\qbasic{4}  &  q\sup3  & 4\\
		&&&&&&&  \\
		& 1   & q\sup2\qbasic{5}   &  q\sup{6}\qcomb{5}{2}  &  q\sup{5}\qcomb{5}{3}  &  q\sup6\qbasic{5}  & q\sup{2} & 5\\	
		&&&&&&&  \\
		1   & q\qbasic{6} & q\sup4\qcomb{6}{2}   &  q\sup{2}\qcomb{6}{3}  &  q\sup{2}\qcomb{6}{4}  &  q\sup{4}\qbasic{6}  & q & 6\\	
		&&&&&&&  \\
		\hline
		6   & 5 & 4  &  3  &  2 &  1  & 0 & l\\	
		\hline
	 \end{array}
	\]}
	\begin{center}
		{\tiny {\bf Figure 1.} Table of the coefficients $\alpha\sub{k,i}$ for $N=7$. Each horizontal row corresponds 		to some
		$0\leq k\leq 6$ and each vertical column corresponds to a fixed $l=(k-i)$.
		The powers of $q$ have been simplified with the identity $q\sup7=1$.}
	\end{center}
	These coefficients can be simplified by using the properties of $q$-numbers. A simple inspection suggests that 
	$\alpha\sub{l}=0$ for $0\leq l\leq N-2$. This is, indeed, the case. Let us write
	{\small\[
		\begin{array}{llr}
			\alpha\sub{l}
			& = \underset{_{l=k-i}}{\overset{_{}}{\sum}}\alpha\sub{k,i}
			   = \underset{_{l=k-i}}{\overset{_{}}{\sum}}\ \ q\sup{i(N-1-k+i)}\cdot\qcomb{k}{i}  & \\
			  && \\
			& = \underset{_{i=0}}{\overset{_{N-l-1}}{\sum}}\ \ q\sup{i(N-l-1)}\cdot\qcomb{l+i}{i} &  \\	
			& & \\		 
			& = \underset{_{i=0}}{\overset{_{s}}{\sum}}\ \ q\sup{is}\cdot\qcomb{N-1-s+i}{i} & \text{take }s=N-l-1 \\	
			& & \\		 
			& = \underset{_{i=0}}{\overset{_{s}}{\sum}}\ \ q\sup{is}\cdot\qcomb{N-1-s+i}{N-1-s}=\beta\sub{s}  & 
			\text{symmetry of combinatorials} \\
		\end{array}
	\]}
	We check that $\beta\sub{s}=\alpha\sub{N-1-s}=0$ for $1\leq s \leq N-1$. For $s=1$,  
	\[
		\beta\sub1=q\sup0+q\sup1\qcomb{N-1}{N-2}=1+q\qbasic{N-1}=\qbasic{N}=0=\alpha\sub{N-2}
	\] 
	Assume that $\beta\sub{s}=0$ for some $s\leq N-2$. Then, 
	{\small\[
			\begin{array}{llr}
				\beta\sub{s+1} & = \underset{_{i=0}}{\overset{_{s+1}}{\sum}}\ \ 
			   		q\sup{i(s+1)}\qcomb{N-2-s+i}{N-2-s} 
					 =\underset{_{i=0}}{\overset{_{s+1}}{\sum}}\ q\sup{is}q\sup{i}\qcomb{N-2-s+i}{N-2-s} 
					 \hskip2mm  & \text{by definition} \\
					&& \\
					& =1+\underset{_{i=1}}{\overset{_{s+1}}{\sum}}\ q\sup{is}
					\left(\qcomb{N-1-s+i}{N-1-s}-\qcomb{N-2-s+i}{N-1-s} \right) 
					 & \text{by \S\ref{lema qcomb}-(3)}\\
					&& \\
					& =1+\underset{_{i=1}}{\overset{_{s+1}}{\sum}}\ q\sup{is}
					\qcomb{N-1-s+i}{N-1-s}-   
					\underset{_{i=1}}{\overset{_{s+1}}{\sum}}\ q\sup{is}
					\qcomb{N-2-s+i}{N-1-s}   & \\
					&& \\
					& =1+\left((-1)+
					     \underset{_{i=0}}{\overset{_{s}}{\sum}}\ q\sup{is}
					\qcomb{N-1-s+i}{N-1-s} + q\sup{(s+1)s}\qcomb{N}{N-1-s}\right) \\
					&& \\   
				  & \hskip1cm - \underset{_{i=1}}{\overset{_{s+1}}{\sum}}\ q\sup{is}
					\qcomb{N-2-s+i}{N-1-s}   & \text{split the first sum} \\
				 && \\
				& =\underset{_{i=0}}{\overset{_{s}}{\sum}}\ q\sup{is}
					\qcomb{N-1-s+i}{N-1-s} 
					-\underset{_{j=0}}{\overset{_{s}}{\sum}}\ q\sup{(j+1)s}
					\qcomb{N-1-s+j}{N-1-s}   & \qbasic{N}=0, i=j+1 \text{(2nd sum)} \\
				 && \\
						& =(1-q\sup{s})\beta\sub{s}=0 & \text{by definition}\\
		\end{array}		
	\]}
	By equation  (\ref{eq homotopy4}), the definition of the Newton's terms \S\ref{Newton's truncation} 
	and a dimension argument on $\sigma$, we deduce that 
	\[
			\underset{_{k=0}}{\overset{_{N-1}}{\sum}}\partial\sup{k}K\partial\sup{N-k-1}(\sigma) 			=\frac{\alpha\sub{N-1}}{\qfac{N-1}}\mc{N}\sup{N-1}(\imath)*\mc{N}\sup{0}(\sigma) =\sigma
	\]
	whenever $n=\dim(\sigma)\geq N-1$, and the whole sum in the left term vanishes when $n<N-1$. In other words,
	\begin{equation}\label{eq homotopy5}
		\underset{_{k=0}}{\overset{_{N-1}}{\sum}}
		\partial\sup{k}K\partial\sup{N-k-1}(\sigma) 			
		=\left\{
		  \begin{array}{lll}
					\sigma & & \dim(\sigma)\geq N-1 \\
					&& \\
					0 && \text{else}	\\	
		\end{array}		
		\right.
	\end{equation}
\edem

\section*{Acknowledgments}
G. Padilla would like to thank Professors E. Becerra, V. Tapia and B. Uribe for some helpful conversations, so as A. Barbosa and D. Maya for their remarks on a previous draft manuscript. This article was partially supported by the Universidad Nacional de Colombia.

\end{document}